\pdfoutput=1

\documentclass[12pt]{article}

\usepackage{amsfonts,amssymb,amsthm,amsmath,latexsym}
\usepackage{subeqnarray,url,cite,bm}
\usepackage{xcolor}

\date{22 March 2024}

\begin{document}

\title{\vspace*{-1cm}
       Continued-fraction characterization of \\
       Stieltjes moment sequences with support in $[\xi, \infty)$
      }

\author{ \\
      \hspace*{-1cm}
      {\large Alan D.~Sokal${}^{1,2}$ and James Walrad${}^1$}
   \\[5mm]
     \hspace*{-1.3cm}
      \normalsize
           ${}^1$Department of Mathematics, University College London,
                    London WC1E 6BT, UK   \\[1mm]
     \hspace*{-2.9cm}
      \normalsize
           ${}^2$Department of Physics, New York University,
                    New York, NY 10003, USA     \\[1mm]
       \\
     \hspace*{-0.5cm}
     {\tt sokal@nyu.edu}, {\tt james.walrad.23@ucl.ac.uk}  \\[1cm]
}

\maketitle
\thispagestyle{empty}   

\vspace*{-5mm}

\begin{abstract}
We give a continued-fraction characterization of Stieltjes moment sequences
for which there exists a representing measure with support in $[\xi, \infty)$.
The proof is elementary.
\end{abstract}


\bigskip
\bigskip
\noindent
{\bf Mathematics Subject Classification (MSC 2020) codes:}
44A60 (Primary);  05A15, 30B70, 30E05 (Secondary).

\clearpage

\newtheorem{theorem}{Theorem}
\newtheorem{proposition}[theorem]{Proposition}
\newtheorem{lemma}[theorem]{Lemma}
\newtheorem{corollary}[theorem]{Corollary}
\newtheorem{definition}[theorem]{Definition}
\newtheorem{conjecture}[theorem]{Conjecture}
\newtheorem{question}[theorem]{Question}
\newtheorem{problem}[theorem]{Problem}
\newtheorem{example}[theorem]{Example}

\renewcommand{\theenumi}{\alph{enumi}}
\renewcommand{\labelenumi}{(\theenumi)}
\def\eop{\hbox{\kern1pt\vrule height6pt width4pt
depth1pt\kern1pt}\medskip}
\def\prf{\par\noindent{\bf Proof.\enspace}\rm}
\def\rmk{\par\medskip\noindent{\bf Remark\enspace}\rm}

\newcommand{\textbfit}[1]{\textbf{\textit{#1}}}

\newcommand{\bigdash}{%
\smallskip\begin{center} \rule{5cm}{0.1mm} \end{center}\smallskip}

\newcommand{\safepar}{ {\protect\hfill\protect\break\hspace*{5mm}} }

\newcommand{\be}{\begin{equation}}
\newcommand{\ee}{\end{equation}}
\newcommand{\<}{\langle}
\renewcommand{\>}{\rangle}
\newcommand{\widebar}{\overline}
\def\reff#1{(\protect\ref{#1})}
\def\spose#1{\hbox to 0pt{#1\hss}}
\def\ltapprox{\mathrel{\spose{\lower 3pt\hbox{$\mathchar"218$}}
    \raise 2.0pt\hbox{$\mathchar"13C$}}}
\def\gtapprox{\mathrel{\spose{\lower 3pt\hbox{$\mathchar"218$}}
    \raise 2.0pt\hbox{$\mathchar"13E$}}}
\def\textprime{${}^\prime$}
\def\proof{\par\medskip\noindent{\sc Proof.\ }}
\def\firstproof{\par\medskip\noindent{\sc First Proof.\ }}
\def\secondproof{\par\medskip\noindent{\sc Second Proof.\ }}
\def\alternateproof{\par\medskip\noindent{\sc Alternate Proof.\ }}
\def\algebraicproof{\par\medskip\noindent{\sc Algebraic Proof.\ }}
\def\combinatorialproof{\par\medskip\noindent{\sc Combinatorial Proof.\ }}
\def\proofof#1{\bigskip\noindent{\sc Proof of #1.\ }}
\def\firstproofof#1{\bigskip\noindent{\sc First Proof of #1.\ }}
\def\secondproofof#1{\bigskip\noindent{\sc Second Proof of #1.\ }}
\def\thirdproofof#1{\bigskip\noindent{\sc Third Proof of #1.\ }}
\def\algebraicproofof#1{\bigskip\noindent{\sc Algebraic Proof of #1.\ }}
\def\combinatorialproofof#1{\bigskip\noindent{\sc Combinatorial Proof of #1.\ }}
\def\sketchofproof{\par\medskip\noindent{\sc Sketch of proof.\ }}
\renewcommand{\qed}{ $\square$ \bigskip}
\newcommand{\myendremark}{ $\blacksquare$ \bigskip}
\def\half{ {1 \over 2} }
\def\third{ {1 \over 3} }
\def\twothird{ {2 \over 3} }
\def\smfrac#1#2{{\textstyle{#1\over #2}}}
\def\smhalf{ {\smfrac{1}{2}} }
\def\smquarter{ {\smfrac{1}{4}} }
\newcommand{\real}{\mathop{\rm Re}\nolimits}
\renewcommand{\Re}{\mathop{\rm Re}\nolimits}
\newcommand{\imag}{\mathop{\rm Im}\nolimits}
\renewcommand{\Im}{\mathop{\rm Im}\nolimits}
\newcommand{\sgn}{\mathop{\rm sgn}\nolimits}
\newcommand{\tr}{\mathop{\rm tr}\nolimits}
\newcommand{\tg}{\mathop{\rm tg}\nolimits}
\newcommand{\supp}{\mathop{\rm supp}\nolimits}
\newcommand{\disc}{\mathop{\rm disc}\nolimits}
\newcommand{\diag}{\mathop{\rm diag}\nolimits}
\newcommand{\csch}{\mathop{\rm csch}\nolimits}
\newcommand{\tridiag}{\mathop{\rm tridiag}\nolimits}
\newcommand{\AZ}{\mathop{\rm AZ}\nolimits}
\newcommand{\perm}{\mathop{\rm perm}\nolimits}
\def\hboxscript#1{ {\hbox{\scriptsize\em #1}} }
\renewcommand{\emptyset}{\varnothing}
\newcommand{\eqdef}{\stackrel{\rm def}{=}}

\newcommand{\restrict}{\upharpoonright}

\newcommand{\compinv}{{\langle -1 \rangle}}   

\newcommand{\scra}{{\mathcal{A}}}
\newcommand{\scrb}{{\mathcal{B}}}
\newcommand{\scrc}{{\mathcal{C}}}
\newcommand{\scrd}{{\mathcal{D}}}
\newcommand{\scre}{{\mathcal{E}}}
\newcommand{\scrf}{{\mathcal{F}}}
\newcommand{\scrg}{{\mathcal{G}}}
\newcommand{\scrh}{{\mathcal{H}}}
\newcommand{\scri}{{\mathcal{I}}}
\newcommand{\scrj}{{\mathcal{J}}}
\newcommand{\scrk}{{\mathcal{K}}}
\newcommand{\scrl}{{\mathcal{L}}}
\newcommand{\scrm}{{\mathcal{M}}}
\newcommand{\scrn}{{\mathcal{N}}}
\newcommand{\scro}{{\mathcal{O}}}
\newcommand{\scrp}{{\mathcal{P}}}
\newcommand{\scrq}{{\mathcal{Q}}}
\newcommand{\scrr}{{\mathcal{R}}}
\newcommand{\scrs}{{\mathcal{S}}}
\newcommand{\scrt}{{\mathcal{T}}}
\newcommand{\scrv}{{\mathcal{V}}}
\newcommand{\scrw}{{\mathcal{W}}}
\newcommand{\scrz}{{\mathcal{Z}}}

\newcommand{\bfa}{{\mathbf{a}}}
\newcommand{\bfb}{{\mathbf{b}}}
\newcommand{\bfc}{{\mathbf{c}}}
\newcommand{\bfd}{{\mathbf{d}}}
\newcommand{\bfe}{{\mathbf{e}}}
\newcommand{\bfj}{{\mathbf{j}}}
\newcommand{\bfi}{{\mathbf{i}}}
\newcommand{\bfk}{{\mathbf{k}}}
\newcommand{\bfl}{{\mathbf{l}}}
\newcommand{\bfm}{{\mathbf{m}}}
\newcommand{\bfx}{{\mathbf{x}}}
\renewcommand{\k}{{\mathbf{k}}}
\newcommand{\n}{{\mathbf{n}}}
\newcommand{\vv}{{\mathbf{v}}}
\newcommand{\bv}{{\mathbf{v}}}
\newcommand{\w}{{\mathbf{w}}}
\newcommand{\x}{{\mathbf{x}}}
\newcommand{\cc}{{\mathbf{c}}}
\newcommand{\zero}{{\mathbf{0}}}
\newcommand{\one}{{\mathbf{1}}}
\newcommand{\bmm}{{\mathbf{m}}}

\newcommand{\ahat}{{\widehat{a}}}
\newcommand{\Zhat}{{\widehat{Z}}}

\newcommand{\C}{{\mathbb C}}
\newcommand{\D}{{\mathbb D}}
\newcommand{\Z}{{\mathbb Z}}
\newcommand{\N}{{\mathbb N}}
\newcommand{\Q}{{\mathbb Q}}
\newcommand{\PP}{{\mathbb P}}
\newcommand{\R}{{\mathbb R}}
\newcommand{\RR}{{\mathbb R}}
\newcommand{\E}{{\mathbb E}}

\newcommand{\Sym}{{\mathfrak{S}}}
\newcommand{\SymB}{{\mathfrak{B}}}
\newcommand{\Alt}{{\mathrm{Alt}}}

\newcommand{\myle}{\preceq}
\newcommand{\myge}{\succeq}
\newcommand{\mygt}{\succ}

\newcommand{\B}{{\sf B}}
\newcommand{\OB}{{\sf OB}}
\newcommand{\OS}{{\sf OS}}
\newcommand{\OO}{{\sf O}}
\newcommand{\SP}{{\sf SP}}
\newcommand{\OSP}{{\sf OSP}}
\newcommand{\Eu}{{\sf Eu}}
\newcommand{\ERR}{{\sf ERR}}
\newcommand{\sfB}{{\sf B}}
\newcommand{\sfD}{{\sf D}}
\newcommand{\sfE}{{\sf E}}
\newcommand{\sfG}{{\sf G}}
\newcommand{\sfJ}{{\sf J}}
\newcommand{\sfP}{{\sf P}}
\newcommand{\sfQ}{{\sf Q}}
\newcommand{\sfS}{{\sf S}}
\newcommand{\sfT}{{\sf T}}
\newcommand{\sfW}{{\sf W}}
\newcommand{\sfMV}{{\sf MV}}
\newcommand{\AMV}{{\sf AMV}}
\newcommand{\BM}{{\sf BM}}
\newcommand{\NC}{{\sf NC}}

\newcommand{\emIB}{{\hbox{\em IB}}}
\newcommand{\emIP}{{\hbox{\em IP}}}
\newcommand{\emOB}{{\hbox{\em OB}}}
\newcommand{\emSC}{{\hbox{\em SC}}}

\newcommand{\stat}{{\rm stat}}
\newcommand{\cs}{{\rm cs}}
\newcommand{\cyc}{{\rm cyc}}
\newcommand{\Asc}{{\rm Asc}}
\newcommand{\asc}{{\rm asc}}
\newcommand{\Des}{{\rm Des}}
\newcommand{\des}{{\rm des}}
\newcommand{\Exc}{{\rm Exc}}
\newcommand{\exc}{{\rm exc}}
\newcommand{\Wex}{{\rm Wex}}
\newcommand{\wex}{{\rm wex}}
\newcommand{\Fix}{{\rm Fix}}
\newcommand{\fix}{{\rm fix}}
\newcommand{\lrmax}{{\rm lrmax}}
\newcommand{\rlmax}{{\rm rlmax}}
\newcommand{\Rec}{{\rm Rec}}
\newcommand{\rec}{{\rm rec}}
\newcommand{\Arec}{{\rm Arec}}
\newcommand{\arec}{{\rm arec}}
\newcommand{\ERec}{{\rm ERec}}
\newcommand{\erec}{{\rm erec}}
\newcommand{\EArec}{{\rm EArec}}
\newcommand{\earec}{{\rm earec}}
\newcommand{\recarec}{{\rm recarec}}
\newcommand{\nonrec}{{\rm nonrec}}
\newcommand{\Cpeak}{{\rm Cpeak}}
\newcommand{\cpeak}{{\rm cpeak}}
\newcommand{\Cval}{{\rm Cval}}
\newcommand{\cval}{{\rm cval}}
\newcommand{\Cdasc}{{\rm Cdasc}}
\newcommand{\cdasc}{{\rm cdasc}}
\newcommand{\Cddes}{{\rm Cddes}}
\newcommand{\cddes}{{\rm cddes}}
\newcommand{\cdrise}{{\rm cdrise}}
\newcommand{\cdfall}{{\rm cdfall}}
\newcommand{\Peak}{{\rm Peak}}
\newcommand{\peak}{{\rm peak}}
\newcommand{\Val}{{\rm Val}}
\newcommand{\val}{{\rm val}}
\newcommand{\Dasc}{{\rm Dasc}}
\newcommand{\dasc}{{\rm dasc}}
\newcommand{\Ddes}{{\rm Ddes}}
\newcommand{\ddes}{{\rm ddes}}
\newcommand{\inv}{{\rm inv}}
\newcommand{\maj}{{\rm maj}}
\newcommand{\rs}{{\rm rs}}
\newcommand{\cross}{{\rm cr}}
\newcommand{\crosshat}{{\widehat{\rm cr}}}
\newcommand{\nest}{{\rm ne}}
\newcommand{\rodd}{{\rm rodd}}
\newcommand{\reven}{{\rm reven}}
\newcommand{\lodd}{{\rm lodd}}
\newcommand{\leven}{{\rm leven}}
\newcommand{\sg}{{\rm sg}}
\newcommand{\bl}{{\rm bl}}
\newcommand{\tran}{{\rm tr}}
\newcommand{\area}{{\rm area}}
\newcommand{\ret}{{\rm ret}}
\newcommand{\peaks}{{\rm peaks}}
\newcommand{\hl}{{\rm hl}}
\newcommand{\sll}{{\rm sl}}
\newcommand{\negg}{{\rm neg}}
\newcommand{\imp}{{\rm imp}}
\newcommand{\osg}{{\rm osg}}
\newcommand{\ons}{{\rm ons}}
\newcommand{\isg}{{\rm isg}}
\newcommand{\ins}{{\rm ins}}
\newcommand{\LL}{{\rm LL}}
\newcommand{\height}{{\rm ht}}
\newcommand{\as}{{\rm as}}

\newcommand{\ba}{{\bm{a}}}
\newcommand{\bahat}{{\widehat{\bm{a}}}}
\newcommand{\sfa}{{{\sf a}}}
\newcommand{\bb}{{\bm{b}}}
\newcommand{\bc}{{\bm{c}}}
\newcommand{\bchat}{{\widehat{\bm{c}}}}
\newcommand{\bd}{{\bm{d}}}
\newcommand{\bee}{{\bm{e}}}
\newcommand{\bff}{{\bm{f}}}
\newcommand{\bg}{{\bm{g}}}
\newcommand{\bh}{{\bm{h}}}
\newcommand{\bll}{{\bm{\ell}}}
\newcommand{\bp}{{\bm{p}}}
\newcommand{\br}{{\bm{r}}}
\newcommand{\bs}{{\bm{s}}}
\newcommand{\bu}{{\bm{u}}}
\newcommand{\bw}{{\bm{w}}}
\newcommand{\bx}{{\bm{x}}}
\newcommand{\by}{{\bm{y}}}
\newcommand{\bz}{{\bm{z}}}
\newcommand{\bA}{{\bm{A}}}
\newcommand{\bB}{{\bm{B}}}
\newcommand{\bC}{{\bm{C}}}
\newcommand{\bE}{{\bm{E}}}
\newcommand{\bF}{{\bm{F}}}
\newcommand{\bG}{{\bm{G}}}
\newcommand{\bH}{{\bm{H}}}
\newcommand{\bI}{{\bm{I}}}
\newcommand{\bJ}{{\bm{J}}}
\newcommand{\bM}{{\bm{M}}}
\newcommand{\bN}{{\bm{N}}}
\newcommand{\bP}{{\bm{P}}}
\newcommand{\bQ}{{\bm{Q}}}
\newcommand{\bS}{{\bm{S}}}
\newcommand{\bT}{{\bm{T}}}
\newcommand{\bW}{{\bm{W}}}
\newcommand{\bX}{{\bm{X}}}
\newcommand{\bIB}{{\bm{IB}}}
\newcommand{\bOB}{{\bm{OB}}}
\newcommand{\bOS}{{\bm{OS}}}
\newcommand{\bERR}{{\bm{ERR}}}
\newcommand{\bSP}{{\bm{SP}}}
\newcommand{\bMV}{{\bm{MV}}}
\newcommand{\bBM}{{\bm{BM}}}
\newcommand{\balpha}{{\bm{\alpha}}}
\newcommand{\bbeta}{{\bm{\beta}}}
\newcommand{\bgamma}{{\bm{\gamma}}}
\newcommand{\bdelta}{{\bm{\delta}}}
\newcommand{\bkappa}{{\bm{\kappa}}}
\newcommand{\bomega}{{\bm{\omega}}}
\newcommand{\bsigma}{{\bm{\sigma}}}
\newcommand{\btau}{{\bm{\tau}}}
\newcommand{\bpsi}{{\bm{\psi}}}
\newcommand{\bzeta}{{\bm{\zeta}}}
\newcommand{\bone}{{\bm{1}}}
\newcommand{\bzero}{{\bm{0}}}

\newcommand{\Cbar}{{\overline{C}}}
\newcommand{\Dbar}{{\overline{D}}}
\newcommand{\dbar}{{\overline{d}}}
\def\Ctilde{{\widetilde{C}}}
\def\Etilde{{\widetilde{E}}}
\def\Ftilde{{\widetilde{F}}}
\def\Gtilde{{\widetilde{G}}}
\def\Htilde{{\widetilde{H}}}
\def\Ptilde{{\widetilde{P}}}
\def\Chat{{\widehat{C}}}
\def\ctilde{{\widetilde{c}}}
\def\zbar{{\overline{Z}}}
\def\pitilde{{\widetilde{\pi}}}

\newcommand{\sech}{{\rm sech}}

%
%
\newcommand{\sn}{{\rm sn}}
\newcommand{\cn}{{\rm cn}}
\newcommand{\dn}{{\rm dn}}
\newcommand{\sm}{{\rm sm}}
\newcommand{\cm}{{\rm cm}}

%
%
\newcommand{\zfz}{ {{}_0 \! F_0} }
\newcommand{\zfo}{ {{}_0  F_1} }
\newcommand{\ofz}{ {{}_1 \! F_0} }
\newcommand{\ofo}{ {{}_1 \! F_1} }
\newcommand{\oft}{ {{}_1 \! F_2} }

%
%
\newcommand{\FHyper}[2]{ {\tensor[_{#1 \!}]{F}{_{#2}}\!} }
\newcommand{\FHYPER}[5]{ {\FHyper{#1}{#2} \!\biggl(
   \!\!\begin{array}{c} #3 \\[1mm] #4 \end{array}\! \bigg|\, #5 \! \biggr)} }
\newcommand{\tfo}{ {\FHyper{2}{1}} }
\newcommand{\tfz}{ {\FHyper{2}{0}} }
\newcommand{\threefz}{ {\FHyper{3}{0}} }
\newcommand{\FHYPERbottomzero}[3]{ {\FHyper{#1}{0} \!\biggl(
   \!\!\begin{array}{c} #2 \\[1mm] \hbox{---} \end{array}\! \bigg|\, #3 \! \biggr)} }
\newcommand{\FHYPERtopzero}[3]{ {\FHyper{0}{#1} \!\biggl(
   \!\!\begin{array}{c} \hbox{---} \\[1mm] #2 \end{array}\! \bigg|\, #3 \! \biggr)} }

\newcommand{\phiHyper}[2]{ {\tensor[_{#1}]{\phi}{_{#2}}} }
\newcommand{\psiHyper}[2]{ {\tensor[_{#1}]{\psi}{_{#2}}} }
\newcommand{\PhiHyper}[2]{ {\tensor[_{#1}]{\Phi}{_{#2}}} }
\newcommand{\PsiHyper}[2]{ {\tensor[_{#1}]{\Psi}{_{#2}}} }
\newcommand{\phiHYPER}[6]{ {\phiHyper{#1}{#2} \!\left(
   \!\!\begin{array}{c} #3 \\ #4 \end{array}\! ;\, #5, \, #6 \! \right)\!} }
\newcommand{\psiHYPER}[6]{ {\psiHyper{#1}{#2} \!\left(
   \!\!\begin{array}{c} #3 \\ #4 \end{array}\! ;\, #5, \, #6 \! \right)} }
\newcommand{\PhiHYPER}[5]{ {\PhiHyper{#1}{#2} \!\left(
   \!\!\begin{array}{c} #3 \\ #4 \end{array}\! ;\, #5 \! \right)\!} }
\newcommand{\PsiHYPER}[5]{ {\PsiHyper{#1}{#2} \!\left(
   \!\!\begin{array}{c} #3 \\ #4 \end{array}\! ;\, #5 \! \right)\!} }
\newcommand{\zerophizero}{ {\phiHyper{0}{0}} }
\newcommand{\ophizero}{ {\phiHyper{1}{0}} }
\newcommand{\zphio}{ {\phiHyper{0}{1}} }
\newcommand{\ophio}{ {\phiHyper{1}{1}} }
\newcommand{\tphio}{ {\phiHyper{2}{1}} }
\newcommand{\tphiz}{ {\phiHyper{2}{0}} }
\newcommand{\tPhio}{ {\PhiHyper{2}{1}} }
\newcommand{\opsio}{ {\psiHyper{1}{1}} }

%
%
\newcommand{\stirlingsubset}[2]{\genfrac{\{}{\}}{0pt}{}{#1}{#2}}
\newcommand{\stirlingcycle}[2]{\genfrac{[}{]}{0pt}{}{#1}{#2}}
\newcommand{\assocstirlingsubset}[3]{{\genfrac{\{}{\}}{0pt}{}{#1}{#2}}_{\! \ge #3}}
\newcommand{\genstirlingsubset}[4]{{\genfrac{\{}{\}}{0pt}{}{#1}{#2}}_{\! #3,#4}}
\newcommand{\euler}[2]{\genfrac{\langle}{\rangle}{0pt}{}{#1}{#2}}
\newcommand{\eulergen}[3]{{\genfrac{\langle}{\rangle}{0pt}{}{#1}{#2}}_{\! #3}}
\newcommand{\eulersecond}[2]{\left\langle\!\! \euler{#1}{#2} \!\!\right\rangle}
\newcommand{\eulersecondgen}[3]{{\left\langle\!\! \euler{#1}{#2} \!\!\right\rangle}_{\! #3}}
\newcommand{\binomvert}[2]{\genfrac{\vert}{\vert}{0pt}{}{#1}{#2}}
\newcommand{\binomsquare}[2]{\genfrac{[}{]}{0pt}{}{#1}{#2}}


\newenvironment{sarray}{
             \textfont0=\scriptfont0
             \scriptfont0=\scriptscriptfont0
             \textfont1=\scriptfont1
             \scriptfont1=\scriptscriptfont1
             \textfont2=\scriptfont2
             \scriptfont2=\scriptscriptfont2
             \textfont3=\scriptfont3
             \scriptfont3=\scriptscriptfont3
           \renewcommand{\arraystretch}{0.7}
           \begin{array}{l}}{\end{array}}

\newenvironment{scarray}{
             \textfont0=\scriptfont0
             \scriptfont0=\scriptscriptfont0
             \textfont1=\scriptfont1
             \scriptfont1=\scriptscriptfont1
             \textfont2=\scriptfont2
             \scriptfont2=\scriptscriptfont2
             \textfont3=\scriptfont3
             \scriptfont3=\scriptscriptfont3
           \renewcommand{\arraystretch}{0.7}
           \begin{array}{c}}{\end{array}}


\newcommand*\circled[1]{\tikz[baseline=(char.base)]{
  \node[shape=circle,draw,inner sep=1pt] (char) {#1};}}
\newcommand{\ostar}{{\circledast}}
\newcommand{\ostarN}{{\,\circledast_{\vphantom{\dot{N}}N}\,}}
\newcommand{\ostarPsi}{{\,\circledast_{\vphantom{\dot{\Psi}}\Psi}\,}}
\newcommand{\starN}{{\,\ast_{\vphantom{\dot{N}}N}\,}}
\newcommand{\starpsi}{{\,\ast_{\vphantom{\dot{\bpsi}}\!\bpsi}\,}}
\newcommand{\starone}{{\,\ast_{\vphantom{\dot{1}}1}\,}}
\newcommand{\startwo}{{\,\ast_{\vphantom{\dot{2}}2}\,}}
\newcommand{\starinfty}{{\,\ast_{\vphantom{\dot{\infty}}\infty}\,}}
\newcommand{\starT}{{\,\ast_{\vphantom{\dot{T}}T}\,}}

\newcommand*{\Perm}[2]{{}^{#1}\!P_{#2}}%

\clearpage

Let us recall that a sequence $\ba = (a_n)_{n \ge 0}$ of real numbers
is called a Hamburger (resp.\ Stieltjes, resp.\ Hausdorff) moment sequence
\cite{Shohat_43,Akhiezer_65,Berg_84,Simon_98,Schmudgen_17}
if there exists a positive measure $\mu$
on $\R$ (resp.\ on $[0,\infty)$, resp.\ on $[0,1]$)
such that $a_n = \int \! x^n \, d\mu(x)$ for all $n \ge 0$;
we call $\mu$ a {\em representing measure}\/ for $\ba$.
A Hamburger (resp.\ Stieltjes) moment sequence is called
H-determinate (resp.\ S-determinate)
if there is a unique such measure $\mu$;
otherwise it is called H-indeterminate (resp.\ S-indeterminate).
Please note that a Stieltjes moment sequence can be
S-determinate but H-indeterminate
\cite[p.~240]{Akhiezer_65} \cite[p.~96]{Simon_98}.

One fundamental characterization of Stieltjes moment sequences
was found by Stieltjes \cite{Stieltjes_1894} in 1894
(see also \cite[pp.~327--329]{Wall_48}):
A sequence $\ba = (a_n)_{n \ge 0}$ of real numbers
is a Stieltjes moment sequence if and only~if
there exist real numbers $c \ge 0$ and $\alpha_1,\alpha_2,\ldots \ge 0$
such that
\be
   \sum_{n=0}^{\infty} a_n t^n
   \;=\;
   \cfrac{c}{1 - \cfrac{\alpha_1 t}{1 - \cfrac{\alpha_2 t}{1 - \cdots}}}
   \label{eq.Sfrac}
\ee
in the sense of formal power series.
(That is, the ordinary generating function
 $f(t) = \sum\limits_{n=0}^{\infty} a_n t^n$
 can be represented as a Stieltjes-type continued fraction
 with nonnegative coefficients.)
Moreover, when $c \neq 0$
the coefficients $\balpha = (\alpha_i)_{i \ge 1}$ are unique
if we make the convention that $\alpha_i = 0$ implies
$\alpha_j = 0$ for all $j > i$;
we shall call such a sequence $\balpha$ {\em standard}\/.

If $K$ is a closed subset of $[0,\infty)$,
we denote by $\scrs_K$ the set of Stieltjes moment sequences
for which there exists a representing measure $\mu$
with support in $K$.\footnote{
   We stress the words ``there exists'':
   it is {\em not}\/ necessarily the case that {\em all}\/
   representing measures have support in $K$, when $K$ is unbounded.
   Indeed, a stunning contrary situation occurs whenever the
   Stieltjes moment sequence is S-indeterminate:
   see Remark~3 at the end of this paper.
}
Can we characterize, in terms of the coefficients $\balpha$,
the sequences in $\scrs_K$?

For the case $K = [0,1]$
--- or by a rescaling, $K = [0,\xi]$ for any $\xi > 0$ ---
this question was answered by Wall \cite[Theorems~4.1 and 6.1]{Wall_40}
in 1940:

\begin{theorem}[Wall 1940]
   \label{thm.wall}
Fix $\xi \ge 0$,
and let $\ba = (a_n)_{n \ge 0}$ be a sequence of real numbers.
Then the following are equivalent:
\begin{itemize}
   \item[(a)] $\ba$ is a Stieltjes moment sequence for which there exists
       a representing measure with support in $[0,\xi]$.
       (In this situation the representing measure is unique.)
   \item[(b)] There exist real numbers $c \ge 0$,
       $g_0 = 0$ and $g_1,g_2,\ldots \in [0,1]$
       such that \reff{eq.Sfrac} holds in the sense of formal power series,
       with
\be
   \alpha_n \;=\;  \xi \, (1 - g_{n-1}) g_n   \;.
\ee
   \item[(c)] Same as (b), but with $g_0 \in [0,1]$.
\end{itemize}
\end{theorem}
Wall \cite{Wall_40,Wall_44} gave two analytic proofs of this result;
and recently one of us \cite{Sokal_wall_hausdorff} gave an elementary proof
of an algebraic/combinatorial flavor.\footnote{
   All of these proofs considered only (a)$\iff$(b).
   But (b)$\implies$(c) is trivially true;
   and (c)$\implies$(a) holds because
   the sequence generated by $g_0 \in (0,1]$
   is bounded above by the corresponding sequence with $g_0 = 0$,
   hence satisfies $a_n \le \xi^n$ according to (b)$\implies$(a).
}

Here we would like to answer the analogous question
for the case $K = [\xi, \infty)$.
We refer to these moment sequences as \textbfit{$\bm{\xi}$-Stieltjes}.
Our result is the following:

\begin{theorem}
   \label{thm.main}
Fix $\xi \ge 0$,
and let $\ba = (a_n)_{n \ge 0}$ be a sequence of real numbers.
Then the following are equivalent:
\begin{itemize}
   \item[(a)] $\ba$ is a $\xi$-Stieltjes moment sequence.
   \item[(b)] There exist real numbers $c \ge 0$, $g_0 = 0$
      and $g_1,g_2,\ldots \ge 0$
such that \reff{eq.Sfrac} holds in the sense of formal power series, with
\begin{subeqnarray}
   \alpha_{2k-1}  & = &  \xi \, (1 + g_{2k-2}) \,+\,  g_{2k-1}  \\[3mm]
   \alpha_{2k}    & = &  {g_{2k-1} g_{2k} \over 1 + g_{2k-2}}
 \label{eq.thm.main}
\end{subeqnarray}
   \item[(c)] Same as (b), but with $g_0 \ge 0$.
\end{itemize}
\end{theorem}


Note in particular that in order for $\ba$
to be a $\xi$-Stieltjes moment sequence,
it is necessary that $\alpha_i \ge \xi$ for all odd $i$;
but it is not sufficient: see Remark~1 below.

Our proof of Theorem~\ref{thm.main}
will use arguments similar to those used in \cite{Sokal_wall_hausdorff}
to prove Theorem~\ref{thm.wall},
but the logic is in many ways simpler
(though the algebra is a bit messier).

Besides Stieltjes-type continued fractions \reff{eq.Sfrac}
[henceforth called S-fractions for short],
we shall also make use of Jacobi-type continued fractions (J-fractions)
\be
   f(t)
   \;=\;
   \cfrac{1}{1 - \gamma_0 t - \cfrac{\beta_1 t^2}{1 - \gamma_1 t - \cfrac{\beta_2 t^2}{1 - \cdots}}}
  \label{eq.Jfrac}
\ee
(always considered as formal power series in the indeterminate $t$).\footnote{
   Our use of the terms ``S-fraction'' and ``J-fraction'' follows
   the general practice in the combinatorial literature,
   starting with Flajolet \cite{Flajolet_80}.
   The classical literature on continued fractions
   \cite{Perron,Wall_48,Jones_80,Lorentzen_92,Cuyt_08}
   generally uses a different terminology.
   For instance, Jones and Thron \cite[pp.~128--129, 386--389]{Jones_80}
   use the term ``regular C-fraction''
   for (a minor variant of) what we~have called an S-fraction,
   and the term ``associated continued fraction''
   for (a minor variant of) what we~have called a J-fraction.
}
We shall need two elementary facts about these continued fractions:

1) The contraction formula:  We have
\be
   \cfrac{1}{1 - \cfrac{\alpha_1 t}{1 - \cfrac{\alpha_2 t}{1 -  \cfrac{\alpha_3 t}{1- \cdots}}}}
   \;\;=\;\;
   \cfrac{1}{1 - \alpha_1 t - \cfrac{\alpha_1 \alpha_2 t^2}{1 - (\alpha_2 + \alpha_3) t - \cfrac{\alpha_3 \alpha_4 t^2}{1 - (\alpha_4 + \alpha_5) t - \cfrac{\alpha_5 \alpha_6 t^2}{1- \cdots}}}}
 \label{eq.contraction_even}
\ee
as an identity between formal power series.
In other words, an S-fraction with coefficients $\balpha$
is equal to a J-fraction with coefficients $\bgamma$ and $\bbeta$, where
\begin{subeqnarray}
   \gamma_0  & = &  \alpha_1
       \slabel{eq.contraction_even.coeffs.a}   \\
   \gamma_n  & = &  \alpha_{2n} + \alpha_{2n+1}  \quad\hbox{for $n \ge 1$}
       \slabel{eq.contraction_even.coeffs.b}   \\
   \beta_n  & = &  \alpha_{2n-1} \, \alpha_{2n}
       \slabel{eq.contraction_even.coeffs.c}
 \label{eq.contraction_even.coeffs}
\end{subeqnarray}
See \cite[pp.~20--22]{Wall_48} for the classic algebraic proof of
the contraction formula \reff{eq.contraction_even};
see \cite[Lemmas~1 and 2]{Dumont_94b}
\cite[proof of Lemma~1]{Dumont_95} \cite[Lemma~4.5]{DiFrancesco_10}
for a very simple variant algebraic proof;
and see \cite[pp.~V-31--V-32]{Viennot_83}
for an enlightening combinatorial proof,
based on Flajolet's \cite{Flajolet_80} combinatorial interpretation
of S-fractions (resp.\ J-fractions)
as generating functions for Dyck (resp.\ Motzkin) paths
with height-dependent weights.

2) Binomial transform:  Fix a real number $\xi$,
and let $\ba = (a_n)_{n \ge 0}$ be a sequence of real numbers.
Then the $\xi$-binomial transform of $\ba$
is defined to be the sequence $\bb = (b_n)_{n \ge 0} \eqdef B_\xi \ba$ given by
\be
   b_n  \;=\; \sum_{k=0}^n \binom{n}{k} \, a_k \, \xi^{n-k}
   \;.
\ee
Note that if $a_n = \int \! x^n \, d\mu(x)$,
then $b_n = \int (x + \xi)^n \, d\mu(x)$.
In other words, if $\ba$ is a Hamburger moment sequence
with representing measure $\mu$,
then $\bb = B_\xi \ba$ is a Hamburger moment sequence
with representing measure $T_\xi \mu$
(the $\xi$-translate of $\mu$).

Now suppose that the ordinary generating function of $\ba$
is given by a J-fraction:
\be
   \sum_{n=0}^{\infty} a_n t^n
   \;=\;
   \cfrac{1}{1 - \gamma_0 t - \cfrac{\beta_1 t^2}{1 - \gamma_1 t - \cfrac{\beta_2 t^2}{1 - \cdots}}}
   \;.
  \label{eq.Jfrac.bis}
\ee
Then the $\xi$-binomial transform $\bb$ of $\ba$ is given by a J-fraction
in which we make the replacement $\gamma_i \to \gamma_i + \xi$:
\be
   \sum_{n=0}^{\infty} b_n t^n
   \;=\;
   \cfrac{1}{1 - (\gamma_0 + \xi) t - \cfrac{\beta_1 t^2}{1 - (\gamma_1 + \xi) t - \cfrac{\beta_2 t^2}{1 - \cdots}}}
   \;.
  \label{eq.Jfrac.bis2}
\ee
See \cite[Proposition~4]{Barry_09} for an algebraic proof
of \reff{eq.Jfrac.bis2};
or see \cite{Sokal_totalpos} for a simple combinatorial proof
based on Flajolet's \cite{Flajolet_80} theory.


\bigskip

\proofof{Theorem~\ref{thm.main}}
When $\xi = 0$ the theorem is trivial,
because any standard nonnegative sequence $\balpha$
can be obtained in the form \reff{eq.thm.main}
with a suitable nonnegative $\bg$.
(We say ``standard'' because, under \reff{eq.thm.main},
 $\alpha_{2k-1} = 0$ forces $\alpha_{2k} = 0$.)
So we assume henceforth that $\xi > 0$.

Note now that a sequence $\ba = (a_n)_{n \ge 0}$
is a $\xi$-Stieltjes moment sequence if and only~if
it is the $\xi$-binomial transform of some Stieltjes moment sequence $\bb$;
and of course this means that $\bb = B_{-\xi} \ba$.

(a)$\implies$(b):
Let $\bb = B_{-\xi} \ba$ be a Stieltjes moment sequence,
and let $\balpha = (\alpha_i)_{i \ge 1}$
be the standard nonnegative coefficients
in the S-fraction \reff{eq.Sfrac} of its ordinary generating function.
Then by the contraction formula \reff{eq.contraction_even},
the ordinary generating function of $\bb$ is also given by a J-fraction
with the coefficients \reff{eq.contraction_even.coeffs}.
Therefore, the ordinary generating function of $\ba = B_\xi \bb$
is given by a J-fraction with the coefficients
\reff{eq.contraction_even.coeffs} replaced by $\gamma_i \to \gamma_i + \xi$:
that is,
\begin{subeqnarray}
   \gamma_0  & = &  \alpha_1 + \xi
       \\
   \gamma_n  & = &  \alpha_{2n} + \alpha_{2n+1} + \xi  \quad\hbox{for $n \ge 1$}
       \\
   \beta_n  & = &  \alpha_{2n-1} \, \alpha_{2n}
 \label{eq.Jfrac.a}
\end{subeqnarray}

Now set $g_0 = 0$,
and define $g_1,g_2,g_3,\ldots \ge 0$ inductively by the rules
\begin{subeqnarray}
   g_{2k-1}  & = &  \alpha_{2k-1}
        \\[2mm]
   g_{2k}  & = &
     \displaystyle
     {1 + g_{2k-2}  \over  \xi \, (1 + g_{2k-2}) \,+\,  g_{2k-1}} \: \alpha_{2k}
 \slabel{def.g.b}
 \label{def.g}
\end{subeqnarray}
Since $\xi > 0$, the prefactor in \reff{def.g.b} is strictly positive;
so the map $\balpha \mapsto \bg$
is a bijection from nonnegative sequences $\balpha$
to nonnegative sequences $\bg$.
(Also $g_i > 0$ if and only~if $\alpha_i > 0$, though we will not need this.)

Now define $\balpha' = (\alpha'_i)_{i \ge 1}$ by \reff{eq.thm.main}, i.e.
\begin{subeqnarray}
   \alpha'_{2k-1}  & = &  \xi \, (1 + g_{2k-2}) \,+\,  g_{2k-1}  \\[3mm]
   \alpha'_{2k}    & = &  {g_{2k-1} g_{2k} \over 1 + g_{2k-2}}
 \label{def.alphaprime}
\end{subeqnarray}
Substituting \reff{def.g} into \reff{def.alphaprime},
we find after some algebra that
\begin{subeqnarray}
   \alpha'_1  & = &  \alpha_1 + \xi  \\
   \alpha'_{2k} + \alpha'_{2k+1}  & = &  \alpha_{2k} + \alpha_{2k+1} + \xi
      \quad\hbox{for $k \ge 1$}
           \\
   \alpha'_{2k-1} \, \alpha'_{2k}  & = &  \alpha_{2k-1} \, \alpha_{2k} 
\end{subeqnarray}
In other words, $\balpha'$ are precisely the S-fraction coefficients
that are related by contraction with the J-fraction coefficients
\reff{eq.Jfrac.a}.
This shows that the ordinary generating function of $\ba$
is given by an S-fraction with the coefficients $\balpha'$.

(b)$\implies$(a):
Given a sequence $\balpha'$ satisfying \reff{def.alphaprime}
for some nonnegative sequence $\bg$,
we can find a nonnegative sequence $\balpha$
satisfying \reff{def.g} because the map $\balpha \mapsto \bg$ is a bijection.
The foregoing argument, read backwards, shows that
the ordinary generating function of $\bb = B_{-\xi} \ba$
has an S-fraction representation \reff{eq.Sfrac}
with nonnegative coefficients $\balpha$,
and hence is a Stieltjes moment sequence.

(b)$\implies$(c) holds trivially.

(c)$\implies$(b):
Let $g_0 > 0$ and $g_1,g_2,\ldots \ge 0$,
and define $\balpha$ by \reff{eq.thm.main}.
We will show that for any $g'_0 \in [0,g_0]$,
there exist $g'_1,g'_2,\ldots \ge 0$
that give rise via \reff{eq.thm.main} to the same coefficients $\balpha$;
in particular, we can choose $g'_0 = 0$. 
Indeed, we will prove inductively that
\begin{subeqnarray}
   0 \;\le\; g'_{2k} & \le &  g_{2k}
      \slabel{eq.induct.2k}  \\[2mm]
   g'_{2k+1} & \ge & g_{2k+1}  \;\ge\; 0
      \slabel{eq.induct.2k+1}
      \label{eq.induct}
\end{subeqnarray}
for all $k \ge 0$.
To see this, observe first that the equations
asserting the equality $\balpha = \balpha'$ are
\begin{subeqnarray}
   \xi g_{2k-2} + g_{2k-1}  & = &  \xi g'_{2k-2} + g'_{2k-1}
           \slabel{eq.ctob.2k-1}  \\[2mm]
   {g_{2k-1} g_{2k} \over 1 + g_{2k-2}}  & = &
      {g'_{2k-1} g'_{2k} \over 1 + g'_{2k-2}}
           \slabel{eq.ctob.2k}
           \label{eq.ctob}
\end{subeqnarray}
for $k \ge 1$.

{\em Base case $k=0$ of \reff{eq.induct}.}
The $k=0$ case of \reff{eq.induct.2k} is the hypothesis $g'_0 \in [0,g_0]$.
And \reff{eq.ctob.2k-1} with $k=1$ gives
\be
   g'_1  \;=\;  g_1 + \xi (g_0 - g'_0)  \;\ge\; g_1
   \;,
\ee
as claimed in the $k=0$ case of \reff{eq.induct.2k+1}.
This proves the base case of the induction.

{\em Inductive step.}
Now assume that (\ref{eq.induct}a,b) hold for a given $k$.
Then using \reff{eq.ctob.2k} with $k \to k+1$, we have
\be
   g'_{2k+2}  \;=\;  g_{2k+2} \: {g_{2k+1} \over g'_{2k+1}}
                              \: {1 + g'_{2k} \over 1 + g_{2k}}
         \;\le\; g_{2k+2}
 \label{eq.induct.2k+2}
\ee
by the inductive hypotheses (\ref{eq.induct}a,b),
provided that $g'_{2k+1} \neq 0$.
But if $g'_{2k+1} = 0$, then the inductive hypothesis \reff{eq.induct.2k+1}
says that $g_{2k+1}  = 0$,
in which case we can satisfy both \reff{eq.ctob.2k} and \reff{eq.induct.2k}
by setting $g'_{2k+2} = 0$.
And using \reff{eq.ctob.2k-1} with $k \to k+2$, we have
\be
   g'_{2k+3}  \;=\;  g_{2k+3} \,+\, \xi (g_{2k+2} - g'_{2k+2})
              \;\ge\; g_{2k+3}
\ee
by \reff{eq.induct.2k+2}.
This completes the inductive step.
\qed

\medskip

{\bf Remarks.}
1. It is easy to see that, although $\alpha_i \ge \xi$ for all odd $i$
is {\em necessary}\/ for $\ba$ to be a $\xi$-Stieltjes moment sequence,
it is not sufficient.
For example, the Catalan numbers
$C_n = \displaystyle \frac{1}{n+1} \binom{2n}{n}$
are a Stieltjes moment sequence with ordinary generating function
\be
   \sum_{n=0}^{\infty} C_n \, t^n
   \;=\;
   \cfrac{1}{1 - \cfrac{t}{1 - \cfrac{t}{1 - \cdots}}}
   \;,
\ee
i.e.~$\alpha_i = 1$ for all $i$.
And the Catalan numbers have a unique representing measure
$(2\pi)^{-1} \,  x^{-1/2} \, (4-x)^{1/2} \, dx$ on the interval $[0,4]$,
so they are not a $\xi$-Stieltjes moment sequence for any $\xi > 0$.
This can also be seen from the criterion of Theorem~\ref{thm.main}:
if $\xi>0$ and $g_0,g_1,g_2,\ldots\ge 0$ satisfy
\begin{subeqnarray}
   1  & = &  \xi \, (1 + g_{2k-2}) \,+\,  g_{2k-1}
      \slabel{eq.catalan.1}  \\[3mm]
   1  & = &  {g_{2k-1} g_{2k} \over 1 + g_{2k-2}}
\end{subeqnarray}
then
$g_{2k} \ge [1- \xi \, (1+ g_{2k-2})] \, g_{2k} = g_{2k-1}g_{2k} = 1+g_{2k-2}$
for all $k$, so $g_{2k} \ge k$.
But for large $k$ this is incompatible with \reff{eq.catalan.1}
whenever $\xi > 0$.

2.  Fix $\xi > 0$, let $\ba$ be a $\xi$-Stieltjes moment sequence,
and let $\balpha$ be the standard nonnegative coefficients
in its S-fraction \reff{eq.Sfrac}.
Now let $I(\ba)$ be the set of $g_0 \ge 0$
for which there exist $g_1, g_2, \ldots \ge 0$
that reproduce the given coefficients $\balpha$ via \reff{eq.thm.main}.
Theorem~\ref{thm.main}(a)$\implies$(b) states that $0 \in I(\ba)$.
Define $g_0^{\rm max} = \sup I(\ba)$.
From eq.~(\ref{eq.thm.main}a) we see that
$g_0^{\rm max} \le \alpha_1/\xi - 1 < \infty$;
and the foregoing proof of (c)$\implies$(b) shows that
$I(\ba)$ is an interval,
either $[0,g_0^{\rm max})$ or $[0,g_0^{\rm max}]$. We can show that the latter is the case: take any sequence $g_0\uparrow g_0^{\rm max}$; then the $g_i'$ are rational functions of $g_0$, which are bounded on both sides by \reff{eq.induct}; therefore they must converge to finite nonnegative limits as $g_0\uparrow g_0^{\rm max}$.

There are obviously cases where $g_0^{\rm max} > 0$;
indeed, we can choose {\em arbitrary}\/ numbers $g_0, g_1, g_2, \ldots \ge 0$
and define $\balpha$ by \reff{eq.thm.main}
and then $\ba$ by \reff{eq.Sfrac};
then clearly $g_0^{\rm max} \ge g_0$.
However, there are also $\xi$-Stieltjes moment sequences $\ba$
for which $g_0^{\rm max} = 0$, as we now show by explicit example.

Consider the sequence $a_n = \smhalf(1^n + 2^n)$,
which is a $1$-Stieltjes moment sequence.
It has the ordinary generating function
\be
   \sum_{n=0}^{\infty} a_n \, t^n
   \;=\;
   {1 \over 2} \, \biggl( {1 \over 1-t} \,+\, {1 \over 1-2t} \biggr)
   \;=\;
   \cfrac{1}{1 - \cfrac{\smfrac{3}{2} t}{1 - \cfrac{\smfrac{1}{6} t}{1 - \smfrac{4}{3} t}}}
   \;,
\ee
so that $\alpha_1 = \smfrac{3}{2}$, $\alpha_2 = \smfrac{1}{6}$,
$\alpha_3 = \smfrac{4}{3}$ and $\alpha_i = 0$ for $i \ge 4$.
Solving \reff{eq.thm.main} with $\xi = 1$ and some fixed $g_0 \ge 0$,
we find
\begin{subeqnarray}
   g_1  & = &  \frac{1}{2} - g_0    \\[1mm]
   g_2  & = &  {1 + g_0 \over 3 (1 - 2g_0)}   \\[1mm]
   g_3  & = &  - \, {g_0 \over 3 (1 - 2g_0)}
\end{subeqnarray}
So $g_1 \ge 0$ requires that $0 \le g_0 \le \half$,
but this is compatible with $g_3 \ge 0$ only if $g_0 = 0$.
This example can of course be rescaled to allow any $\xi > 0$.

3. Let $\ba$ be any Stieltjes moment sequence that is
S-indeterminate, i.e.~for which there exists more than one
representing measure supported on $[0,\infty)$.
Then it is a remarkable and surprising (at least to us) fact that
$\ba$ is {\em always}\/ a $\xi$-Stieltjes moment sequence for some $\xi > 0$.
This follows from a result of Berg and Valent
\cite[Lemma~2.2.1 and Remark~2.2.2]{Berg_94b},
as follows:

Let $\ba$ be any Stieltjes moment sequence that is
H-indeterminate,
and let $A,B,C,D$ be the corresponding Nevanlinna functions
\cite[eqns.~(2.10)/(2.12)]{Berg_94b}.
Then Berg and Valent \cite[Lemma~2.2.1]{Berg_94b} show that
\be
   \alpha  \,\eqdef\,  \lim_{x \to -\infty} {D(x) \over B(x)}
   \;\hbox{ exists and is } \le 0   \;.
\ee
They furthermore observe \cite[Remark~2.2.2]{Berg_94b}
that the N-extremal measure $\nu_t$ is supported on $[0,\infty)$
if and only if $t \in [\alpha,0]$,
and that $\alpha < 0$ if and only~if $\ba$ is S-indeterminate.
(As always,
all the N-extremal measures $(\nu_t)_{t \in \R \cup \{\infty\}}$
have disjoint supports.)
Furthermore, when $\alpha < 0$ it follows from \cite[Lemma~2.2.1]{Berg_94b}
that $\xi(t) \eqdef \inf\supp\nu_t$
is nonnegative, continuous, real-analytic and strictly decreasing on
the interval $[\alpha,0]$, with $\xi(0) = 0$;
in particular, we can achieve a value $\xi^* \eqdef \xi(\alpha) > 0$.
Further reasoning \cite{Berg_private}
based on \cite{Pedersen_95,Pedersen_97}
shows that $\nu_\alpha$ is the {\em unique}\/ representing measure $\mu$
having $\inf\supp\mu = \xi^*$,
and that there does not exist any representing measure
having $\inf\supp\mu > \xi^*$.

It is an interesting problem, for concrete examples of
S-indeterminate moment sequences $\ba$,
to determine explicitly the maximum ``spectral gap'' $\xi^*$.
For Stieltjes' \cite[section~56]{Stieltjes_1894} famous example
$a_n = q^{-n(n-1)/2}$ with $q \in (0,1)$,
this computation was in essence performed
by Christiansen \cite{Christiansen_03}.
\myendremark

\section*{Acknowledgments}

We are extremely grateful to Christian Berg for discussions
concerning Remark~3.

\end{document}